\documentclass[%
 aip,
cp,  
 amsmath,amssymb,
 reprint,%
]{revtex4-2}

\usepackage{eufrak}
\usepackage{graphicx}
\usepackage{dcolumn}
\usepackage{bm}
\usepackage[hidelinks]{hyperref}
\usepackage[capitalize,nameinlink]{cleveref}[0.19]
\usepackage{xcolor}
\usepackage{soul}

\usepackage[utf8]{inputenc}
\usepackage[T1]{fontenc}
\usepackage{mathptmx}

\newcommand{\boltz}{\textup{k}_{\textup{B}}}

\newcommand{\lambdaSq}{\chi}

\newcommand{\tim}{t}

\newcommand{\StateVec}{\boldsymbol{w}}
\newcommand{\SVi}[1][i]{\StateVec_{#1}}

\newcommand{\SVspecies}{\SVi[\species]}

\newcommand{\dwspecies}{\partial_{\SVspecies}}
\newcommand{\flx}{\boldsymbol{f}}
\newcommand{\flxi}[1][i]{\flx_{#1}}

\newcommand{\flxspecies}{\flxi[\species]}
\newcommand{\jcb}{A}
\newcommand{\jcbi}[1][i]{\jcb_{#1}}

\newcommand{\jcbspecies}{\jcbi[\species]}

\newcommand{\St}{\boldsymbol{S}}
\newcommand{\Sti}[1][i]{\St_{#1}}

\newcommand{\Stspecies}{\Sti[\species]}
\newcommand{\DscrtFlx}{F}

\newcommand{\bDscrtFlx}{\boldsymbol{F}}

\newcommand{\AprxFlx}{\mathcal{F}}
\newcommand{\AprxFlxRus}{\mathcal{F}^{Rus}}
\newcommand{\AprxFlxHLL}{\mathcal{F}^{HLL}}

\newcommand{\PoissonOp}{\mathcal{P}_\lambdaSq}
\newcommand{\IzOp}{\mathcal{I}_{0}}
\newcommand{\IzOpone}{\mathcal{I}_{1}}

\newcommand{\Mscalingi}[1][i]{M^{s}_{#1}}

\newcommand{\Mscalingion}{\Mscalingi[\ion]}

\newcommand{\Domain}{\Omega}

\newcommand{\eigvalspecies}{\lambda^{\species}}
\newcommand{\eigvalelec}{\lambda^{\elec}}

\newcommand{\eigvalspeciesi}[1][i]{\eigvalspecies_{#1}}
\newcommand{\eigvaleleci}[1][i]{\eigvalelec_{#1}}

\newcommand{\eigvalspeciespm}{\eigvalspeciesi[\pm]}
\newcommand{\lbdmaxspecies}{\eigvalspeciesi[max]}
\newcommand{\lbdmaxelec}{\eigvaleleci[max]}

\newcommand{\lbdmaxspeciesjph}{\eigvalspeciesi[{max_{j+1/2}}]}

\newcommand{\lbdmaxelecnjph}{\lambda^{\elec,n}_{max_{j+1/2}}}
\newcommand{\lbdmaxelecnjmh}{\lambda^{\elec,n}_{max_{j-1/2}}}

\newcommand{\HLLbpspecies}{b^{\species}_{+}}
\newcommand{\HLLbmspecies}{b^{\species}_{-}}

\newcommand{\HLLbpspeciesnjph}{b^{\species,n}_{+,j+1/2}}
\newcommand{\HLLbmspeciesnjph}{b^{\species,n}_{-,j+1/2}}

\newcommand{\HLLbpion}{b^{\ion}_{+}}
\newcommand{\HLLbmion}{b^{\ion}_{-}}
\newcommand{\HLLbpmion}{b^{\ion}_{\pm}}

\newcommand{\SVspeciesSj}{\StateVec^{*}_{\species_{j}}}

\newcommand{\SVspeciesnj}{\StateVec^{n}_{\species_{j}}}
\newcommand{\SVspeciesnjp}{\StateVec^{n}_{\species_{j+1}}}

\newcommand{\SVspeciesnNp}{\StateVec^{n}_{\species_{N+1}}}
\newcommand{\SVspeciesno}{\StateVec^{n}_{\species_{0}}}

\newcommand{\SVionnNp}{\StateVec^{n}_{\ion_{N+1}}}
\newcommand{\SVionno}{\StateVec^{n}_{\ion_{0}}}
\newcommand{\SVionnN}{\StateVec^{n}_{\ion_{N}}}
\newcommand{\SVionnone}{\StateVec^{n}_{\ion_{1}}}

\newcommand{\SVspeciesnonej}{\StateVec^{n_1}_{\species_{j}}}
\newcommand{\SVspeciesntwoj}{\StateVec^{n_2}_{\species_{j}}}
\newcommand{\SVspeciesnthreej}{\StateVec^{n_3}_{\species_{j}}}
\newcommand{\SVspeciesnfourj}{\StateVec^{n_4}_{\species_{j}}}

\newcommand{\SVspeciesnpj}{\StateVec^{n+1}_{\species_{j}}}

\newcommand{\SVspeciesR}{\StateVec_{\species_R}}
\newcommand{\SVspeciesL}{\StateVec_{\species_L}}

\newcommand{\DFlxspeciesnjph}{\DscrtFlx^{n}_{\species_{j+1/2}}}
\newcommand{\DFlxspeciesnjmh}{\DscrtFlx^{n}_{\species_{j-1/2}}}

\newcommand{\DFlxionno}{\DscrtFlx^{n}_{\ion_{0}}}

\newcommand{\DFlxionnone}{\DscrtFlx^{n}_{\ion_{1}}}

\newcommand{\DFlxionnN}{\DscrtFlx^{n}_{\ion_{N}}}

\newcommand{\DFlxionnNp}{\DscrtFlx^{n}_{\ion_{N+1}}}

\newcommand{\bMFlxspeciesntwojph}{\bDscrtFlx^{MUSCL,n_2}_{\species_{j+1/2}}}
\newcommand{\bMFlxspeciesntwojmh}{\bDscrtFlx^{MUSCL,n_2}_{\species_{j-1/2}}}

\newcommand{\phibnj}{\bar{\phi}^{n}_{j}}
\newcommand{\phibnpj}{\bar{\phi}^{n+1}_{j}}
\newcommand{\phibnpjp}{\bar{\phi}^{n+1}_{j+1}}
\newcommand{\phibnpjm}{\bar{\phi}^{n+1}_{j-1}}
\newcommand{\phibnjpm}{\bar{\phi}^{n}_{j\pm1}}
\newcommand{\phibnpjpm}{\bar{\phi}^{n+1}_{j\pm1}}
\newcommand{\phibnpo}{\bar{\phi}^{n+1}_{0}}
\newcommand{\phibnpone}{\bar{\phi}^{n+1}_{1}}
\newcommand{\phibnpN}{\bar{\phi}^{n+1}_{N}}
\newcommand{\phibnpNp}{\bar{\phi}^{n+1}_{N+1}}
\newcommand{\phibVecnj}{\left( \phibnj \right)_{1\leq j\leq N}}
\newcommand{\phibVecnpj}{\left( \phibnpj \right)_{1\leq j\leq N}}

\newcommand{\nelecbnj}{\bar{n}^{n}_{\elec_j}}
\newcommand{\nelecbnonej}{\bar{n}^{n_1}_{\elec_j}}

\newcommand{\nelecbnthreej}{\bar{n}^{n_3}_{\elec_j}}
\newcommand{\nelecbnfourj}{\bar{n}^{n_4}_{\elec_j}}
\newcommand{\nelecbnjp}{\bar{n}^{n}_{\elec_{j+1}}}
\newcommand{\nelecbnjm}{\bar{n}^{n}_{\elec_{j-1}}}
\newcommand{\nionbnj}{\bar{n}^{n}_{\ion_j}}

\newcommand{\nspeciesbnj}{\bar{n}^{n}_{\species_j}}

\newcommand{\nelecbnpj}{\bar{n}^{n+1}_{\elec_j}}
\newcommand{\nionbnpj}{\bar{n}^{n+1}_{\ion_j}}
\newcommand{\nspeciesbnpj}{\bar{n}^{n+1}_{\species_j}}
\newcommand{\nelecbSj}{\bar{n}^{*}_{\elec_j}}
\newcommand{\nionbSj}{\bar{n}^{*}_{\ion_j}}
\newcommand{\nspeciesbsSj}{\bar{n}^{**}_{\species_j}}
\newcommand{\nelecbsSj}{\bar{n}^{**}_{\elec_j}}

\newcommand{\nelecVecbSj}{\left(\bar{n}^{*}_{\elec_j}\right)_{1\leq j\leq N}}
\newcommand{\nelecVecbnj}{\left(\bar{n}^{n}_{\elec_j}\right)_{1\leq j\leq N}}
\newcommand{\nelecVecbnonej}{\left(\bar{n}^{n_1}_{\elec_j}\right)_{1\leq j\leq N}}

\newcommand{\nelecVecbnthreej}{\left(\bar{n}^{n_3}_{\elec_j}\right)_{1\leq j\leq N}}
\newcommand{\nelecVecbnfourj}{\left(\bar{n}^{n_4}_{\elec_j}\right)_{1\leq j\leq N}}

\newcommand{\nelecnionVecbSj}{\left(\bar{n}^{*}_{\elec_j},\bar{n}^{*}_{\ion_j}\right)_{1\leq j\leq N}}
\newcommand{\nelecnionVecbnthreej}{\left(\bar{n}^{n_3}_{\elec_j},\bar{n}^{n_3}_{\ion_j}\right)_{1\leq j\leq N}}

\newcommand{\nelecnionVecbSnj}{\left(\bar{n}^{n}_{\ion_j},\bar{n}^{*}_{\elec_j},\bar{n}^{*}_{\ion_j}\right)_{1\leq j\leq N}}
\newcommand{\nelecnionVecbsSnj}{\left(\bar{n}^{n}_{\ion_j},\bar{n}^{**}_{\elec_j},\bar{n}^{**}_{\ion_j}\right)_{1\leq j\leq N}}

\newcommand{\DamkS}{\bar{\nu}^{iz,*}}
\newcommand{\DamksS}{\bar{\nu}^{iz,**}}
\newcommand{\Damkn}{\bar{\nu}^{iz,n}}
\newcommand{\Damknone}{\bar{\nu}^{iz,n_1}}

\newcommand{\Damknthree}{\bar{\nu}^{iz,n_3}}
\newcommand{\Damknfour}{\bar{\nu}^{iz,n_4}}

\newcommand{\nuionbSj}{\bar{\nu}^{*}_{\ion_j}}

\newcommand{\nuspeciesbSj}{\bar{\nu}^{*}_{\species_j}}

\newcommand{\nUebnjp}{\left(\bar{n}_\elec\bar{u}_\elec\right)^{n}_{j+1}}
\newcommand{\nUebnjm}{\left(\bar{n}_\elec\bar{u}_\elec\right)^{n}_{j-1}}

\newcommand{\nUspeciesbnpj}{\left(\bar{n}_\species\bar{u}_\species\right)^{n+1}_{j}}

\newcommand{\nUspeciesbnpjph}{\left(\bar{n}_\species\bar{u}_\species\right)^{n+1}_{j+1/2}}
\newcommand{\nUspeciesbnpjmh}{\left(\bar{n}_\species\bar{u}_\species\right)^{n+1}_{j-1/2}}
\newcommand{\nUspeciesbSj}{\left(\bar{n}_\species\bar{u}_\species\right)^{*}_{j}}

\newcommand{\nUebno}{\left(\bar{n}_\elec\bar{u}_\elec\right)_{0}}
\newcommand{\nUebnone}{\left(\bar{n}_\elec\bar{u}_\elec\right)_{1}}
\newcommand{\nUebnN}{\left(\bar{n}_\elec\bar{u}_\elec\right)_{N}}
\newcommand{\nUebnNp}{\left(\bar{n}_\elec\bar{u}_\elec\right)_{N+1}}

\newcommand{\nUibnone}{\left(\bar{n}_\ion\bar{u}_\ion\right)^{n}_{1}}
\newcommand{\nUibnN}{\left(\bar{n}_\ion\bar{u}_\ion\right)^{n}_{N}}

\newcommand{\nUiVecbnoneN}{\left(\left(\bar{n}_\ion\bar{u}_\ion\right)^{n}_{j}\right)_{j\in\{1,N\}}}
\newcommand{\nUiVecbnonejoneN}{\left(\left(\bar{n}_\ion\bar{u}_\ion\right)^{n_1}_{j}\right)_{j\in\{1,N\}}}

\newcommand{\nUiVecbnthreejoneN}{\left(\left(\bar{n}_\ion\bar{u}_\ion\right)^{n_3}_{j}\right)_{j\in\{1,N\}}}
\newcommand{\nUiVecbnfourjoneN}{\left(\left(\bar{n}_\ion\bar{u}_\ion\right)^{n_4}_{j}\right)_{j\in\{1,N\}}}

\newcommand{\nelecGhost}{\bar{n}_{\elec_G}}
\newcommand{\nelecBndry}{\bar{n}_{\elec_B}}
\newcommand{\nelecbno}{\bar{n}_{\elec_0}}
\newcommand{\nelecbnone}{\bar{n}_{\elec_1}}
\newcommand{\nelecbnN}{\bar{n}_{\elec_N}}
\newcommand{\nelecbnNp}{\bar{n}_{\elec_{N+1}}}
\newcommand{\uelecGhost}{\bar{u}_{\elec_G}}
\newcommand{\uelecBndry}{\bar{u}_{\elec_B}}
\newcommand{\phiGhost}{\bar{\phi}_{G}}
\newcommand{\phiBndry}{\bar{\phi}_{B}}
\newcommand{\OutNorm}{e_n}

\newcommand{\uwall}{u_{wall}}
\newcommand{\uwallb}{\bar{u}_{wall}}

\newcommand{\Dt}{\Delta t}
\newcommand{\DtHLL}{\Delta t_{HLL}}
\newcommand{\DtRus}{\Delta t_{Rus}}
\newcommand{\DtConv}{\Delta t_{Conv}}
\newcommand{\DtSt}{\Delta t_{St}}
\newcommand{\DtIz}{\Delta t_{Iz}}
\newcommand{\DtCol}{\Delta t_{Col}}
\newcommand{\DtPoisson}{\Delta t_{\phi}}
\newcommand{\Dx}{\Delta x}
\newcommand{\DxPoisson}{\Delta x_\phi}
\newcommand{\xprime}{x^{\prime}}

\newcommand{\dt}{\partial_{\tim}}

\newcommand{\dx}{\partial_{x}}
\newcommand{\dxx}{\partial^2_{xx}}

\newcommand{\dtb}{\partial_{\bar{\tim}}}

\newcommand{\dxb}{\partial_{\bar{x}}}
\newcommand{\dxxb}{\partial^2_{\bar{x}\bar{x}}}

\newcommand{\nablab}{\bar{\nabla}}

\newcommand{\xjph}{x_{j+1/2}}
\newcommand{\xjmh}{x_{j-1/2}}

\newcommand{\ee}{e}
\newcommand{\species}{\alpha}
\newcommand{\elec}{\EuFrak{e}}
\newcommand{\neut}{\EuFrak{n}}

\newcommand{\ion}{\EuFrak{i}}

\newcommand{\Damk}{\bar{\nu}^{iz}}

\newcommand{\OmPe}{\omega_{p\elec}}
\newcommand{\OmPeb}{\bar{\omega}_{p\elec}}

\newcommand{\debye}{\lambda_D}
\newcommand{\debyeb}{\bar{\lambda}_D}
\newcommand{\Vsheath}{V_f}
\newcommand{\Vpresheath}{V_s}
\newcommand{\Vsheathb}{\bar{V}_f}
\newcommand{\Vpresheathb}{\bar{V}_s}

\newcommand{\Lo}{L_0}

\newcommand{\lambdai}[1][i]{\lambda_{#1}}

\newcommand{\lambdaion}{\lambdai[\ion]}

\newcommand{\mi}[1][i]{m_{#1}}
\newcommand{\me}{\mi[\elec]}
\newcommand{\mion}{\mi[\ion]}
\newcommand{\mspecies}{\mi[\species]}

\newcommand{\symboln}{n}
\newcommand{\nii}[1][i]{\symboln_{#1}}
\newcommand{\nelec}{\nii[\elec]}

\newcommand{\nspecies}{\nii[\species]}
\newcommand{\nion}{\nii[\ion]}

\newcommand{\no}{\nii[0]}

\newcommand{\nui}[1][i]{\nu_{#1}}
\newcommand{\nue}{\nui[\elec]}
\newcommand{\nuion}{\nui[\ion]}
\newcommand{\nuspecies}{\nui[\species]}
\newcommand{\nuib}[1][i]{\bar{\nu}_{#1}}
\newcommand{\nueb}{\nuib[\elec]}
\newcommand{\nuionb}{\nuib[\ion]}
\newcommand{\nuspeciesb}{\nuib[\species]}
\newcommand{\nuionization}{\nu^{iz}}

\newcommand{\pii}[1][i]{p_{#1}}
\newcommand{\pelec}{\pii[\elec]}
\newcommand{\pspecies}{\pii[\species]}

\newcommand{\pion}{\pii[\ion]}
\newcommand{\phio}{\phi_0}
\newcommand{\phib}{\bar{\phi}}
\newcommand{\permvac}{\epsilon_0}

\newcommand{\qi}[1][i]{q_{#1}}
\newcommand{\qelec}{\qi[\elec]}
\newcommand{\qion}{\qi[\ion]}
\newcommand{\qspecies}{\qi[\species]}

\newcommand{\qib}[1][i]{\bar{q}_{#1}}
\newcommand{\qelecb}{\qib[\elec]}
\newcommand{\qionb}{\qib[\ion]}
\newcommand{\qspeciesb}{\qib[\species]}

\newcommand{\rhoib}[1][i]{\bar{n}_{#1}}
\newcommand{\rhoeb}{\rhoib[\elec]}
\newcommand{\rhoionb}{\rhoib[\ion]}
\newcommand{\rhospeciesb}{\rhoib[\species]}

\newcommand{\timeo}{t_0}
\newcommand{\timei}{t_i}
\newcommand{\timef}{t_f}
\newcommand{\timen}{t^{n}}
\newcommand{\timenp}{t^{n+1}}
\newcommand{\Ti}[1][i]{T_{#1}}
\newcommand{\Te}{\Ti[\elec]}

\newcommand{\Tion}{\Ti[\ion]}

\newcommand{\Tspecies}{\Ti[\species]}
\newcommand{\Tpo}{\Ti[0]}

\newcommand{\vthi}[1][i]{c_{#1}}
\newcommand{\vthelec}{\vthi[\elec]}
\newcommand{\vthion}{\vthi[\ion]}
\newcommand{\vthspecies}{\vthi[\species]}
\newcommand{\vthb}[1][i]{\bar{c}_{#1}}
\newcommand{\vtheb}{\vthb[\elec]}
\newcommand{\vthib}{\vthb[\ion]}
\newcommand{\vthibeff}{\vthb[\ion_{eff}]}
\newcommand{\vthspeciesb}{\vthb[\species]}

\newcommand{\uo}{u_0}

\newcommand{\uBohmb}{\bar{u}_B}
\newcommand{\ui}[1][i]{u_{#1}}
\newcommand{\ue}{\ui[\elec]}
\newcommand{\uion}{\ui[\ion]}
\newcommand{\uspecies}{\ui[\species]}

\newcommand{\uib}[1][i]{\bar{u}_{#1}}
\newcommand{\ueb}{\uib[\elec]}

\newcommand{\uRoespeciesb}{\uib[Roe]^{\species}}

\newcommand{\uRoeionb}{\uib[Roe]^{\ion}}
\newcommand{\uspeciesLb}{\uib[\species_L]}
\newcommand{\uspeciesRb}{\uib[\species_R]}
\newcommand{\uelecLb}{\uib[\elec_L]}
\newcommand{\uelecRb}{\uib[\elec_R]}

\newcommand{\uionb}{\uib[\ion]}

\newcommand{\uionbSj}{\uib[\ion_{j}]^{*}}
\newcommand{\uspeciesb}{\uib[\species]}

\newcommand{\sigmane}{\sigma_{\elec \neut}}
\newcommand{\sigmanion}{\sigma_{\ion \neut}}

\newcommand{\rhospeciesbR}{\rhoib[{\species_{R}}]}
\newcommand{\rhospeciesbL}{\rhoib[{\species_{L}}]}

\begin{document}

\title{Numerical challenges in the simulation of 1D bounded low-temperature plasmas with charge separation in various collisional regimes}

\author{Louis Reboul} 
 \email[Corresponding author: ]{louis.reboul@polytechnique.edu}
\author{Marc Massot}
  \email{marc.massot@polytechnique.edu}
\affiliation{
    CMAP, CNRS, \'Ecole polytechnique, Institut Polytechnique de Paris, Route de Saclay, 91128 Palaiseau Cedex, France
}
\author{Alejandro Alvarez Laguna}%
 \email{alejandro.alvarez-laguna@lpp.polytechnique.fr}
\affiliation{
  Laboratoire de Physique des Plasmas (LPP), CNRS, Sorbonne Université,  \'Ecole Polytechnique, Institut Polytechnique de Paris, 91120 Palaiseau, France.
}

\date{30 November 2022} 

\begin{abstract}
We study a 1D geometry of a plasma confined between two conducting floating walls with applications to laboratory plasmas. These plasmas are characterized by a quasi-neutral bulk that is joined to the wall by a thin boundary layer called sheath that is positively charged. Although analytical solutions are available in the sheath and the pre-sheath, joining the two areas by one analytical solution is still an open problem which requires the numerical resolution of the fluid equations coupled to Poisson equation. Current numerical schemes use high-order discretizations to correctly capture the electron current in the sheath, presenting unsatisfactory results in the boundary layer and they are not adapted to all the possible collisional regimes. In this work, we identify the main numerical challenges that arise when attempting the simulations of such configuration and we propose explanations for the observed phenomena via numerical analysis. We propose a numerical scheme with controlled diffusion as well as new discrete boundary conditions that address the identified issues.
\end{abstract}

\maketitle

Submitted to: API Conference Proceedings.

\section{Introduction}
We consider a 1D plasma bounded between two infinity conductive plates. The density profile in such configuration is well-known in the literature \cite{Chabert,Lieberman} consisting of a quasi-neutral bulk of plasma in the center, and a thin charged region in the boundary layers close to the walls, called sheaths. There is no analytical solution for this problem that matches continuously both the sheath and the pre-sheath, resulting in a need for numerical approximate solutions. The stationary equations of such a configuration are not well posed and consequently a time resolved method must be employed. In the present work, we consider the isothermal fluid equations for ions and electrons coupled to Poisson equation for the electrostatic potential. Previous work \cite{crispel2007asymptotic, laguna2020asymptotic, alvarez2020plasma} have found that classical numerical schemes, wether at first or second order accuracy in time and space, do not yield satisfactory numerical solutions, especially for the profile of electron current (or equivalently electron momentum). 
The objective of the present work is the following: (1) demonstrate that classical numerical methods fail to achieve reasonable accuracy on this set of equations, (2) to analyze the reasons of this lack of accuracies and stability and (3) to propose a simple variation of classical numerical schemes to obtain enhanced accuracy. Finally, the relevance of our new approach is assessed by a number of numerical experiments under different collisional regimes, gas composition and space discretizations.

\section*{Model}
In order to model the 1D plasma-sheath transition, we use a multi-fluid model that considers ions and electrons as two separate fluids interacting via the Lorentz force. Such model can be obtained by taking the two first moment of the kinetic equations \cite{braginskii1965transport}. In this work, we use the isothermal assumption to close the equations. The resulting isothermal Euler-Poisson equations, written here respectively in dimensional (left) and non-dimensional (right) form, read:
\begin{subequations}\label{eq:ScldEulerPoisson} 
  \begin{align}
      \dt \nelec + \dx \left(\nelec \ue\right) &= \nuionization \nelec,
      &\dtb \rhoeb + \dxb \left(\rhoeb \ueb\right) &= \Damk \rhoeb\label{eq:ScldMassElec},\\
      \me\dt \left(\nelec \ue\right) + \dx \left(\me\nelec \ue^2+\pelec\right)&=-\qelec \nelec \dx \phi-\me\nue\nelec\ue,
      &\dtb \left(\rhoeb \ueb\right) + \dxb \left(\rhoeb\left( \ueb^2+\varepsilon^{-1}\right)\right)&=\frac{\rhoeb \dxb \phib}{\varepsilon}-\nueb\rhoeb  \ueb\label{eq:ScldMomElec},\\
      \dt \nion+ \dx \left(\nion \ui\right) &= \nuionization \nelec,
      &\dtb \rhoib+ \dxb \left(\rhoib \uib\right) &= \Damk \rhoeb\label{eq:ScldMassIon},\\
      \mion\dt \left(\nion \ui\right) + \dx \left(\mion\nion \ui^2+\pion\right)&=-\qion \nion \dx \phi-\mion\nuion\nion\ui,
      &\dtb \left(\rhoib\uib\right) + \dxb \left(\rhoib\left( \uib^2+\kappa^{\frac{1}{2}}\right)\right)&=-\rhoib \dxb \phib-\nuionb\rhoionb  \uionb\label{eq:ScldMomIon},\\
      -\dxx \phi &= \frac{\qelec \nelec + \qion \nion}{\epsilon_0},
      &\dxxb \phib &= \chi^{-1} \left(\rhoeb - \rhoib\right)\label{eq:ScldPoisson}.
  \end{align}
\end{subequations}
In the following, we will use the notation $\species\in\left\{\elec,\ion\right\}$ to index quantities related to a species, with $\ion$ denoting ion quantities and $\elec$ electron quantities. On the left-hand side of \cref{eq:ScldEulerPoisson} we have the dimensional equations, with $\nspecies$ referring to number densities,  $\uspecies$ to the macroscopic velocities, $\mspecies$ to the masses of particles, $\pspecies$ to pressure, $\qspecies$ to the charges (here $\qelec = +\ee$ and $\qion = -\ee$), $\phi$ to the electric potential, $\nuspecies$ to the collision rates, $\nuionization$ to the ionization rate (the computation of these various rates is detailed below) and $\permvac$ to vacuum permittivity. Pressure terms are computed assuming the law of perfect gases $\pelec = \nelec \boltz \Te$, $\pion  = \nion  \boltz \Tion$, where $\boltz$ is the Boltzmann constant and $\Tspecies$ are the temperature, assumed to be constant in space and time.

On the right hand side of \cref{eq:ScldEulerPoisson} we have normalized the equations using the Bohm velocity, $\uo = \sqrt{\boltz\Te/\mion}$, that is the velocity at which ions enter the sheath. This choice is adequate for our setting since its purpose is precisely to include the boundary layer that forms at the walls, namely the sheath. Regarding the other characteristic quantities required to scale the equations, we choose $\no$ defined as the average density of ions, $\Tpo = \Te$ the temperature of electrons and $\Lo$ the distance between the two conducting walls.
From these values, other relevant reference quantities can be derived, such as the reference time $\timeo = \Lo/\uo$ and the reference electric potential $\phio = \boltz \Te / e$. Then setting $\rhospeciesb = \nspecies/\no$, $\uspeciesb = \uspecies/\uo$, $\phib = \phi/\phio$, $\Damk = \timeo \nuionization$, $\nueb = \timeo \nue$, $\nuionb = \timeo \nuion$ as well as $\dtb = \timeo \dt$ and $\dxb = \Lo \dx$, one can obtain the normalized equations (right of \cref{eq:ScldEulerPoisson}) from the dimensional equations (left of \cref{eq:ScldEulerPoisson}).

Three small parameters have been defined during this operation: the mass ratio $\varepsilon = \me/\mion$, the temperature ratio $\kappa = \Tion/\Te$ and $\lambdaSq = \left(\boltz \Te\permvac\right) /\left(e^2  \no\Lo^2\right) = \debye^2/\Lo^2$ the square of the normalized Debye length. These parameters define different asymptotic regimes of \cref{eq:ScldEulerPoisson}, as presented in \cite{laguna2020asymptotic}. The limit $\varepsilon\to0$, represents the limit where the inertia of electrons is negligible. In this limit, one can formally equate the terms factored by $\varepsilon^{-1}$ in \cref{eq:ScldMomElec} to obtain the well-known Boltzmann relation $\dxb \rhoeb = \rhoeb \dxb \phib$. Alternatively, the limit $\chi\to0$, represents the limit where the Debye length is negligible as compared to the size of the problem. This limit corresponds to the classic quasi-neutral regime where $\rhoeb \approx \rhoionb$. However, these asymptotic regimes do not correspond to the plasma inside the sheath and, for this reason, all the scales need to remain in the model to capture the plasma-sheath transition. 


Regarding the boundary conditions, we assume that electrons are absorbed by the conducting walls with no secondary electron emission. We assume that electron at the wall have a truncated Maxwellian velocity distribution function, which is compatible with the fact that no electrons are emitted from the wall.  From this velocity distribution one can obtain that the outgoing thermal flux at the boundaries of the domain is $\left| \rhoeb \ueb \right| = \rhoeb/\sqrt{2\pi\varepsilon}$. This boundary condition is usually imposed to study the sheaths with fluid models, as explained in \cite{Chabert}. Imposing these boundary conditions also justifies the assertion that electron inertia in the sheaths cannot be neglected, since at the walls we have $\ueb \propto 1/\sqrt{\varepsilon}$. Boundary conditions are also required to have a well-posed Poisson equation (\cref*{eq:ScldPoisson}), and we choose to set the walls at electric potential zero in this configuration.

The computation of the ionization rate is derived under the assumption that ionization compensates for the loss of ions at the wall, in order to have a steady-state solution. This results in a constant number of ions in the domain, i.e., $\dtb \int_{\Domain} \rhoib\text{d}\bar{V} = 0$, with $\Domain = [0,1]$ being the normalized domain. This condition is used to obtained the ionization frequency by integrating spatially \cref*{eq:ScldMassIon}:
\begin{equation}\label{eq:IzRate}
  \int_{\Domain} \nablab.\left(\rhoib\uionb\right)\text{d}\bar{V} = \int_{\Domain} \rhoeb \Damk\text{d}\bar{V}
  \Leftrightarrow \Damk= \frac{\int_{\partial \bar{\Omega} }\rhoib\uionb\cdot\text{d}\bar{S}}{\int_{\Domain}\rhoeb\text{d}\bar{V}}.
\end{equation}

Finally the elastic collisions rates between the charged species (ions and electrons) and the neutral gas are computed assuming constant collision cross sections and Maxwellian velocity distributions for both ions and electrons, following the computations from \cite{Chabert}:
\begin{equation}\label{eq:colls}
  \nuionb = 
\frac{\Lo}{\lambdaion}\sqrt{\frac{8\kappa}{\pi}+\frac{\pi^2}{4}\uion^2},
\quad
\nueb=\frac{\sigmane}{\sigmanion} \frac{\Lo}{\lambdaion}\sqrt{\frac{8}{\pi \varepsilon}}.
\end{equation}
Note that in this work, we neglect the Coulomb collisions, which is justified in weakly-ionized plasmas where $n_{\ion}\sim n_\elec \ll n_\neut$, where $n_\neut$ is the neutral species density. In \cref{{eq:colls}}, we have used the ion mean free path that is defined as $\lambdaion = (n_\neut  \sigmanion)^{-1}$, where $\sigmanion$ is the ion-neutral momentum-transfer cross section. In the next section we present classical approaches to solve numerically \cref*{eq:ScldEulerPoisson}.


\section{Some classically used numerical schemes}
We rewrite \cref*{eq:ScldEulerPoisson} under the form:
\begin{equation}\label{eq:HypSt}
  \dtb \SVspecies+ \dxb \flxspecies\left(\SVspecies\right) = \Stspecies\left(\SVspecies\right),
\end{equation}
with $\species \in \left\{ \elec, \ion \right\}$ denoting the species, ions and electrons, and
\begin{equation}\label{eq:SVflxStspecies}
  \SVspecies =
  \begin{pmatrix}
    \rhospeciesb\\
    \rhospeciesb\uspeciesb
  \end{pmatrix}, \text{ }
  \flxspecies\left( \SVspecies \right) = 
  \begin{pmatrix}
    \rhospeciesb\uspeciesb\\
    \rhospeciesb\uspeciesb^{2}+\rhospeciesb\vthspeciesb^{2}
  \end{pmatrix}, \text{ }
  \Stspecies \left( \SVspecies \right) = 
  \begin{pmatrix}
    \Damk\rhoeb\\
    -\qspeciesb \dxb \phib-\nuspeciesb\rhospeciesb\uspeciesb
  \end{pmatrix}
\end{equation}
with the speed of sound $\vthspecies = \sqrt{\boltz\Tspecies/\mspecies}$, that in normalized form reads $\vtheb = \varepsilon^{-1/2}$ and $\vthib = \kappa^{1/2}$ and with $\qelecb = -1$ and $\qionb = +1$. Written in this form, \cref*{eq:HypSt} fall under the framework of hyperbolic systems of conservation laws with source terms. The Jacobian matrix $\jcbspecies = \dwspecies \flxspecies$ for the flux has two real eigenvalues $\eigvalspeciespm = \uspeciesb\pm\vthspeciesb$. Therefore, if we set the source term at zero, $\Stspecies = 0$, we retrieve a hyperbolic system, see for instance \cite{godlewski2013numerical}, namely in this case the isothermal Euler equations. Hyperbolic systems can develop singularities, i.e. shocks, in finite time, a technical difficulty that must be properly addressed by the numerical method, otherwise spurious oscillations can appear around discontinuities. Finite volume methods \cite{toro2013riemann,leveque2002finite} are classically employed to overcome these aspects, and consequently we will use this approach for the left-hand side of \cref*{eq:HypSt}.

The numerical scheme approximates the solution $\SVspecies$ on a mesh. We denote this approximation $\SVspeciesnj$ of the quantity $\int_{C_j} \SVspecies\left(x,y,t^n\right)\textup{d}x$ with $C_j = [\xjmh,\xjph]$, $\xjph=j\Dx$, $1\leq j \leq N$, $N\geq 1$, and $t^n = n\Dt$, $n\geq 1$. A fully-explicit scheme is not allowed for stability reasons (see \cite{fabre1992stability}), and the electric potential $\phi$ must be computed in an implicit manner.

A possible approach to achieve stability is to use a first-order Lie operator splitting (see \cite{leveque2002finite,duarte:tel-00667857,DUARTE2015129,duarte:hal-00573043}):
\begin{subequations}\label{eq:Lie1splitting}
  \begin{align}
    &\text{{\bfseries Step 1: }Convective step: }&\SVspeciesSj&=\SVspeciesnj-\frac{\Dt}{\Dx}\left(\DFlxspeciesnjph - \DFlxspeciesnjmh\right),\label{eq:Lie1ConvStep}\\
    &\text{{\bfseries Step 2: }Computation of Poisson equation: }&\phibVecnpj &= \PoissonOp \left(\nelecnionVecbSj\right),\label{eq:Lie1PoissonEq}\\
    &\text{{\bfseries Step 3: }Computation of the ionization rate frequency:}&\DamkS&=\IzOp\left( \nelecnionVecbSnj \right),\label{eq:Lie1IzComp}\\
    &\text{{\bfseries Step 4: }Integration of the source terms: }& \SVspeciesnpj &= \SVspeciesSj +\Dt \Stspecies \left( \SVspeciesSj, \DamkS, \nelecbSj, \phibnpjpm \right),\label{eq:Lie1St}
  \end{align}
\end{subequations}
where $\PoissonOp$ (resp. $\IzOp$) is a discrete operator that yields an approximate solution to the Poisson equation (resp. computes the ionization rate).

We choose this type of splitting because it is known to be quite robust \cite{toro2013riemann,leveque2002finite,hundsdorfer2003numerical} since it decouples the hyperbolic and source term parts and therefore stability can be established by taking the minimal time step between the one obtained via the CFL condition of the convective step and the one used for the source term step seen here as an ODE problem. We now detail each of these steps.

For the convective step (\cref*{eq:Lie1ConvStep}), one must give the approximate value of the fluxes of the form $\DFlxspeciesnjph = \AprxFlx\left( \SVspeciesnj, \SVspeciesnjp \right)$ to be used at the interfaces. In this work we consider two approximate Riemann solvers:
\begin{subequations}\label{eq:ApproxRiemanSolvers}
  \begin{align}
    &\text{Rusanov solver:}&\AprxFlxRus\left( \SVspeciesL, \SVspeciesR \right) &= \frac{\flxspecies\left( \SVspeciesR \right) + \flxspecies\left( \SVspeciesL \right)}{2}-\frac{\lbdmaxspecies}{2}\left( \SVspeciesR - \SVspeciesL \right),\label{eq:RusRS}\\
    &\text{HLL solver:}&\AprxFlxHLL\left( \SVspeciesL, \SVspeciesR \right)&= \frac{\HLLbpspecies\flxspecies\left( \SVspeciesL \right)-\HLLbmspecies\flxspecies\left( \SVspeciesR \right)}{\HLLbpspecies-\HLLbmspecies}+\frac{\HLLbpspecies\HLLbmspecies}{\HLLbpspecies-\HLLbmspecies}\left( \SVspeciesR - \SVspeciesL \right),
  \end{align}
\end{subequations}
where $\lbdmaxspecies = \max\left( \left|\uspeciesLb\right|, \left|\uspeciesRb\right| \right)+\vthspeciesb$, that is the modulus of the maximum eigenvalue of $\jcbspecies$, and $\HLLbpspecies = \max\left( 0, \uRoespeciesb +\vthspeciesb \right)$ and $\HLLbmspecies = \min\left( 0, \uRoespeciesb -\vthspeciesb \right)$ with $\uRoespeciesb = \left( \sqrt{\rhospeciesbL}\uspeciesLb+\sqrt{\rhospeciesbR}\uspeciesRb \right)/\left( \sqrt{\rhospeciesbL}+\sqrt{\rhospeciesbR} \right)$. We will compare the performances of these Riemann solvers in the next section.

A key issue is related to the treatment of boundary conditions. One can note that when $j=1$ or $j=N$ the convective step \cref*{eq:Lie1ConvStep} needs approximate the values $\SVspeciesno$ and $\SVspeciesnNp$ which are outside of the domain. We need to use the boundary conditions of our continuous problem \cref*{eq:ScldEulerPoisson} to fill these cells at index $j=0$ and $j=N+1$, which are commonly called ghost cells \cite{toro2013riemann,leveque2002finite,hundsdorfer2003numerical}.

In the case of low-temperature bounded plasmas we follow the classical approach developed in \cite{alvarez2020plasma}. For ions, we use the fact that in practice $\kappa\ll1$, i.e. $\vthion\ll \uo$, to assume that they are in a supersonic regime when they reach the boundary of the domain within the sheath and therefore we can take an upwind approximation of the form $\SVionno = \SVionnone$ and $\SVionnNp = \SVionnN$. For electrons, the situation is the opposite, as we have $\uwallb := 1/\sqrt{2\pi\varepsilon} < \varepsilon^{-1/2} = \vtheb$, electron will be in a subsonic, or low-Mach, regime, even in the sheath. Consequently, information propagates in both directions and it is not possible to extrapolate both the density and the momentum from within the domain. We choose to extrapolate the electron density, on which we have no clear boundary condition in the continuous problem, from the domain, that is we take $\nelecbno = \nelecbnone$ and $\nelecbnNp = \nelecbnN$. For the momentum, we use the boundary condition at the wall $\left| \rhoeb \ueb \right| = \rhoeb/\sqrt{2\pi\varepsilon}$, which leads to:
\begin{equation*}
  \frac{\nUebno+\nUebnone}{2} = -\nelecbnone \uwallb, \quad \frac{\nUebnN+\nUebnNp}{2} = +\nelecbnN \uwallb,
\end{equation*}
and complete the convective step.

To obtain an approximate solution of the Poisson equation we use the discrete operator $\PoissonOp$ that yields the solution of the linear system defined, for all $j\in\left\{1,...,N\right\}$, by the relation:
\begin{equation}\label{eq:DiscretePoisson}
  \frac{\phibnpjp-2\phibnpj+\phibnpjm}{\Dx^{2}} = \frac{\nelecbSj-\nionbSj}{\lambdaSq}.
\end{equation}
As in the convective step, for $j=1,N$ we need to set values for $\phibnpo$ and $\phibnpNp$. We use once more the boundary conditions of the continuous problem \cref*{eq:ScldPoisson}, leading to:
\begin{equation*}
  \frac{\phibnpo+\phibnpone}{2} = 0, \quad \frac{\phibnpN+\phibnpNp}{2}=0.
\end{equation*}
We use the upper index $n+1$ instead of $*$ because complete scheme leads to the equality $\nelecbnpj-\nionbnpj = \nelecbSj-\nionbSj$ as we will show later in this section.

For the computation of the ionization rate, we use a formula that ensures that the loss of ion mass at the walls is compensated by ionization \cite{alvarez2020plasma}, leading to the following formula for the ionization operator $\IzOp$:
\begin{equation}\label{eq:CompIz}
  \DamkS = \frac{\sum_j \nionbnj-\sum_j \nionbSj}{\Dt \sum_j \nelecbSj}
\end{equation}
Because the convective step is conservative and because of the boundary conditions that we have chosen for ions this formula is strictly equivalent to:
\begin{equation}\label{eq:IsOpone}
  \DamkS = \IzOpone\left( \nUiVecbnoneN, \nelecVecbSj \right) :=\frac{\nUibnN-\nUibnone}{\sum_j \nelecbSj \Dx}
\end{equation}
which is a natural discrete form of the continuous \cref*{eq:IzRate}.

The collision rate for electrons is constant in time and space and for the ion collision rate we simply use formula  (\ref{eq:colls}):
\begin{equation}\label{eq:CompColIon}
  \nuionbSj = \frac{\Lo}{\lambdaion}\sqrt{\frac{8\kappa}{\pi}+\frac{\pi^2}{4}\left( \uionbSj \right)^2}
\end{equation}

Finally, for the source term step, we use the formula:
\begin{equation}
  \Stspecies \left( \SVspeciesSj, \DamkS, \nelecbSj, \phibnpjpm \right) = 
  \begin{pmatrix}
    \DamkS \nelecbSj\\
    -\qspecies\frac{\phibnpjp-\phibnpjm}{2\Dx}-\nuspeciesbSj\nUspeciesbSj
  \end{pmatrix}.
\end{equation}
 We can now see that $\nelecbnpj-\nionbnpj = \nelecbSj-\nionbSj$ which confirms that the electric potential is indeed computed in an implicit manner, as required in \cite{fabre1992stability} for numerical stability when coupling isothermal Euler equation to Poisson equation.

Before moving on to simulations with this numerical scheme, we present its stability constraints. First for the convective step, the CFL condition to ensure that local Riemann problems at the interfaces do not interact during the time step is for each solver (see for instance \cite{toro2013riemann,leveque2002finite}):
\begin{equation}\label{cnd:ConvStep}
  \Dt \leq \DtRus :=\frac{\Dx}{2\max_{j}\lbdmaxspeciesjph}, \quad \Dt \leq \DtHLL := \frac{\Dx}{2\max_{j}\left( \max\left( \HLLbmspeciesnjph,\HLLbpspeciesnjph \right) \right)}
\end{equation}
For the source term, the stability conditions
\begin{subequations}\label{cnd:St}
  \begin{align}
    &\text{Poisson stability condition:}& \Dx \leq& \DxPoisson := \sqrt{\lambdaSq},\label{cnd:StPoisson}\\
    &\text{Ionization stability condition:}& \Dt \leq& \DtIz := 1/\DamkS,\label{cnd:StIz}\\
    &\text{Collision stabtility condition}& \Dt \leq& \DtCol := 1/\max\left( \nueb,\nuionbSj \right).\label{cnd:StCol}
  \end{align}
\end{subequations}
One can note than combining conditions from \cref*{cnd:ConvStep} and \cref*{cnd:StPoisson} we obtain:
\begin{equation}\label{cnd:PlasmaFreq}
  \Dt \leq \DtPoisson := \sqrt{\varepsilon\lambdaSq} = \OmPeb^{-1}
\end{equation}
where $\OmPeb = \OmPe \timeo$ is the normalized electron plasma frequency with $\OmPe = \vthelec/\debye = \sqrt{e^{2}\no/\left(\me\permvac\right)}$, that is the characteristic frequency at which electrons oscillate within the plasma due to electrostatic waves, satisfying the stability condition from \cite{fabre1992stability}.

In order to obtain stability, the time step simply needs to satisfy $\Dt \leq \min\left( \DtConv,\DtSt \right)$ with the convective time step $\DtConv = \DtHLL$ or $\DtConv = \DtRus$, depending on the approximate Riemann solver that is used, and the source term time step $\DtSt = \min(\DtIz,\DtCol, \DtPoisson)$.

We have considered so far a rather very classical approach, for which the literature seems to indicate the stable and robust behavior expected for a first order method.
\section*{Accuracy loss and stability issues of the chosen classical approach}
The purpose of this part is to identify, from a set of well-chosen numerical experiments on a given set-up, a list of three key issues pertaining accuracy and stability the scheme described in the previous section is suffering from.
\subsection*{A well chosen configuration for several numerical experiments}
In all simulations of this paper, we set $\no = 10^{15}$ $m^{-3}$, $\Tpo = 1.2\times 10^{5}$ K (typical low-temperature low-pressure plasma parameters, see \cite{Chabert}), leading to $\debye \approx 8\times 10^{-4}$ m and $\phio = 10$ V. As a reference, in non-collisional test cases, we can calculate the theoretical electric potential drop in the sheath $\Vsheath$ and in the pre-sheath $\Vpresheath$, see \cite{Chabert}. These magnitudes read:
\begin{equation}\label{eq:Vdrop}
  \Vsheath = \frac{1}{2}\frac{\boltz\Te}{\ee}\ln\left( \frac{2\pi \me}{\mion} \right), \quad \Vpresheath = \frac{1}{2}\frac{\boltz\Te}{\ee}.
\end{equation}
The non-dimensional form of these potential drops reads $\Vsheathb = \ln\left( 2\pi\varepsilon \right)/2$ and $\Vpresheathb = 1/2$.

The first test case we consider is collisionless and has the following parameters: $\varepsilon = 1/1836$ (simplified hydrogen gas, $\uo \approx 3 $ m.s$^{-1}$), $\kappa = 0.0025$ ($\Ti = 300$ K), $\debyeb=0.02$ ($\Lo \approx 4$ cm and $\timeo \approx 1.2~\mu$s). We choose a gas with a small atomic mass to show that the inaccuracies that we now demonstrate are present with virtually any type of gas (as hydrogen is the lightest atom). We choose an electron-ion temperature ratio typical of low-temperature plasmas. As a first step, we consider collisionless plasmas (no collisions with the neutral gas). We will consider the impact on the accuracy of the collisional terms later in this paper. The simulation is initialized at time $\timei = 0$ with constant in space number densities $\rhoeb = \rhoionb = 1$ and null velocity $\ueb = \uionb = 0$ and the final time of simulation is set at $\timef = 4\timeo$, for which a steady state is reached.

\subsection*{Identification of several problems from numerical experiments}
\begin{figure}[htbp]
  \centering
  \includegraphics[width=1.\columnwidth]{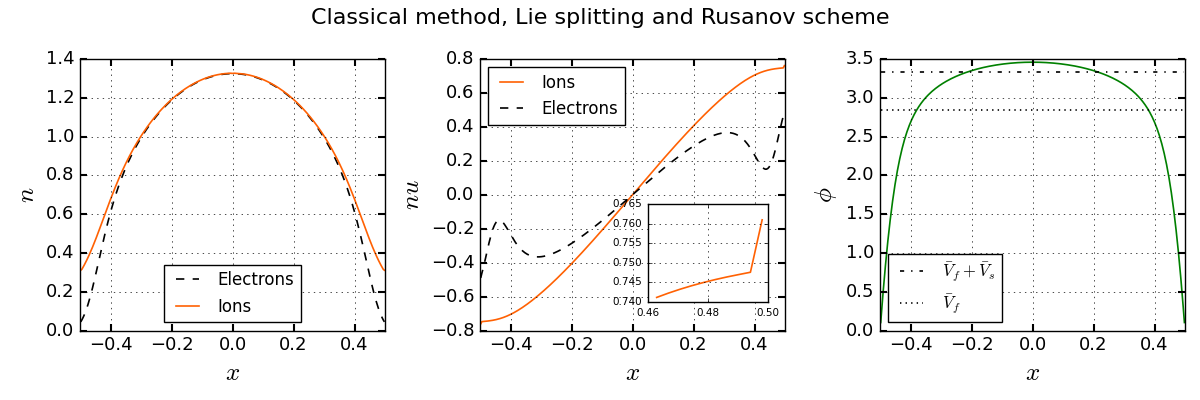}
  \caption{Number density, momentum and electric potential profiles after steady state is reached. Plasma parameters: $\varepsilon = 1/1836$, $\kappa = 0.0025$, $\debyeb = 0.02$. Numerical scheme: Lie splitting, Rusanov solver for both ions and electrons. Simulation parameters: $N = 256$ cells, $\Dx \approx 4\times 10^{-3}$, $\Dt \approx 2\times 10^{-5}$.}
  \label{fig:1stTCe1Rusi1RusLie1}
\end{figure}

First, we use the classical first-order scheme with the Rusanov approximate solver for both ions and electrons, i.e., $\DFlxspeciesnjph = \AprxFlxRus\left( \SVspeciesnj, \SVspeciesnjp \right)$. It is assumed from the literature \cite{toro2013riemann,leveque2002finite,hundsdorfer2003numerical} that this approximate solver is the most robust. The result of the simulation with $N=256$ cells, are displayed in \cref*{fig:1stTCe1Rusi1RusLie1}.
As it can be seen, the number density profiles have the expected bell shape, with a quasi-neutral bulk in the center and charged region at the walls in the sheath, that extends over a few Debye lengths. However, the profile of electron flux seems to be very inaccurate (far superior to the expected first-order precision with $\Dx \approx 4\times 10^{-3}$). As it can be inferred from the model equations, at steady state $\dtb \rhoeb = \dtb \rhoionb = 0$ and consequently by subtracting \cref*{eq:ScldMomIon} to \cref*{eq:ScldMomElec} one obtains $\dxb \left( \rhoeb\ueb - \rhoionb\uib \right) = 0$. Because of the symmetry of the considered configuration we have $\ueb = \uib = 0$ at $x = 0$ and as a result $\rhoeb \ueb  = \rhoionb \uionb$ in the entire domain. In other words, the method should yield an ambipolar profile (i.e., both electron and ion fluxes are the same). This is clearly not the case in \cref*{fig:1stTCe1Rusi1RusLie1}, and the first objective of this work is both to analyze and to remove this large inaccuracies. The electric potential profile does not match accurately the expected potential drop in the sheath and the pre-sheath, which is a possible consequence of the aforementioned inaccuracies.

Another unexpected behavior of the numerical solution is that the profile of ion flux presents a sharp gradient in the last cell close to the wall, as displayed on the zoomed section on the center graph of \cref*{fig:1stTCe1Rusi1RusLie1}. This undesired effect is due to the fact that the boundary conditions that are currently employed implement the value at the ghost cell as $\DFlxionno = \DFlxionnone$ and $\DFlxionnNp = \DFlxionnN$. This is equivalent to suddenly switching from Rusanov to HLL only on the exterior faces of the two cells at the boundary of the domain (index $j=1$ and $j=N$). Not using the same approximate Riemann solver for each face of these cells touching the boundaries is likely the cause of the observed jump in ion momentum.

\begin{figure}[htbp]
  \centering
  \includegraphics[width=1.\columnwidth]{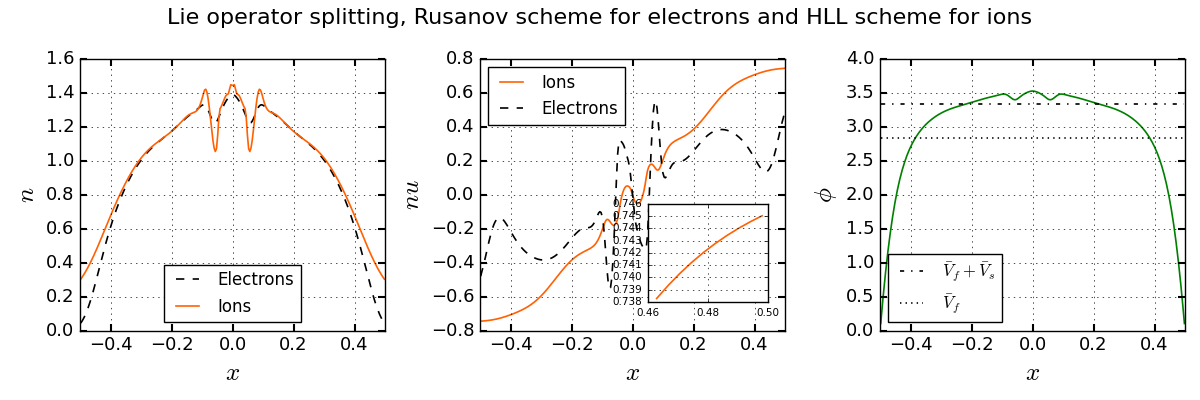}
  \caption{Number density, momentum and electric potential profiles at final time $\timef = 4\timeo$. Plasma parameters: $\varepsilon = 1/1836$, $\kappa = 0.0025$, $\debyeb = 0.02$. Numerical scheme: Lie splitting, Rusanov solver for electron and HLL solver for ions and electrons. Simulation parameters: $N = 256$ cells, $\Dx \approx 4\times 10^{-3}$, $\Dt \approx 2\times 10^{-5}$.}
  \label{fig:1stTCe1Rusi1HLLLie1}
\end{figure}

Switching to the HLL scheme for electron everywhere in the domain is useless since we observed similar inaccuracies in the sheath as with the previous scheme (figures are not presented for the sake of conciseness). Besides, switching to HLL for ions, which could seem natural to resolve the problem at the boundaries, leads to an unstable behavior, as displayed in \cref*{fig:1stTCe1Rusi1HLLLie1} (we have kept the Rusanov scheme for electrons to isolate the source of the unstable behavior).
As a result, and although the profile of ion momentum is indeed smoother in the boundary cells, it is not possible to use the HLL solver for the ion fluid.
The numerical analysis of this phenomenon will be provided at the end of this section. In the next section we will propose an approach allowing to fix this issue.

\begin{figure}[htbp]
  \centering
  \includegraphics[width=1.\columnwidth]{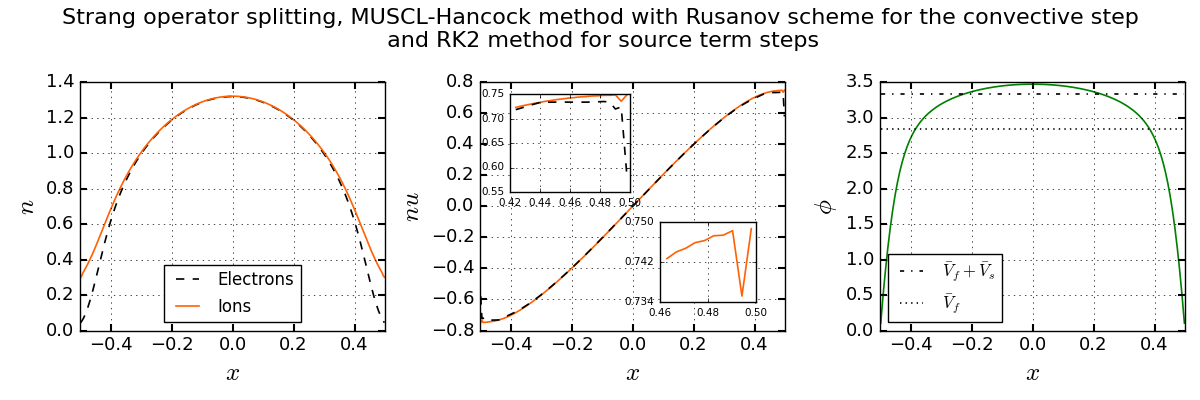}
  \caption{Number density, momentum and electric potential profiles at final time $\timef = 4\timeo$. Plasma parameters: $\varepsilon = 1/1836$, $\kappa = 0.0025$, $\debyeb = 0.02$. Numerical scheme: Strang splitting, RK2 for the source term steps, MUSCL-Hancock methods with minmod slope limiters for the convective step and Rusanov solver for both electrons and ions. Simulation parameters: $N = 256$ cells, $\Dx \approx 4\times 10^{-3}$, $\Dt = 2\times 10^{-5}$.}
  \label{fig:1stTCe2Rusi2RusStrang}
\end{figure}

An alternative method that is aimed to reduce the inaccuracies observed in the sheath area in \cref*{fig:1stTCe1Rusi1RusLie1} is to increase the precision of the scheme and for example to switch to a second-order discretization. We consider a Strang operator splitting, with Runge-Kutta second order method for the source term steps and MUSCL-Hancock method with minmod slope limiters for the convective step, as described in the appendix. The results obtained using the Rusanov solver for both ions and electrons are displayed in \cref*{fig:1stTCe2Rusi2RusStrang}.
We can first observe that there is a much better match between the momentum of ions and electrons. However, when we zoom in the sheath area (top left-hand side of the center graph of \cref*{fig:1stTCe2Rusi2RusStrang}), we can see that there is still a mismatch between the ion and electron fluxes. This error becomes larger for gases with heavier ions (e.g. helium or argon), corresponding to a lower value of $\varepsilon$. Another problem appears in the electron flux that presents a sharp gradient on the wall. Understanding the source of this boundary problem as well as proposing new discrete boundary conditions that remove this problem is the third objective of this paper. We can also note that although the solution appears to be more accurate in the domain of simulation, there is clearly a mismatch between the maximum of the electric potential that is computed by the method and the expected theoretical value.

A final observation from the zoomed section at the bottom right-hand side of the center graph of \cref*{fig:1stTCe2Rusi2RusStrang} is that not only is there a sharp variation for ion momentum at boundary because of the Rusanov scheme, but we can also see small oscillations in the sheath area. We believe they correspond to ion acoustic instabilities that are artificially triggered by the numerical scheme, and originate in the fact that the second-order method that we have chosen does not entirely comply with the stability requirements established in \cite{fabre1992stability}.
We point out that all classical manner to achieve second order in time and space suffer from the same flaws, which is that the electric potential must at some point be computed at a time that is in between $\timen$ and $\timenp$ and therefore is not fully implicit. Consequently, similar oscillations were observed, for low values of $\kappa$, with other second order methods such as the method of lines. This problem of reaching second order in time and space will not be tackled in the present work.

To conclude this section, we have identified three challenges that we wish to analyse and address:
\begin{enumerate}
  \item There are strong inaccuracies in the electron flux profile in the whole domain and mainly in the sheath, both for the first-order method and to a lesser extent for the second-order method.
  \item The HLL approximate solver is unstable when used for ions at low values of $\kappa$.
  \item At the boundary of the domain, the discrete boundary conditions for electrons do not seem to be consistent with the problem.
\end{enumerate}

\subsection*{Numerical analysis of the identified numerical problems of the classical schemes}
\subsubsection*{Electron flux error in the sheath}
We start by analyzing problem number 1, that is the inaccuracies of the methods in the sheath area. At steady state, we have $\nelecbnpj = \nelecbnj$, so using the Rusanov scheme yields 
\begin{equation}\label{eq:NelecDstst}
  \frac{\nUebnjp-\nUebnjm}{2\Dx}-\Dx\frac{1}{\Dx}\left( \lbdmaxelecnjph \frac{\nelecbnjp-\nelecbnj}{\Dx}-\lbdmaxelecnjmh\frac{\nelecbnj-\nelecbnjm}{\Dx} \right) = \DamkS \nelecbSj.
\end{equation}
Using the theory of modified equations (see for instance \cite{leveque2002finite}), which essentially means injecting a continuous solution into \cref*{eq:NelecDstst} and conducting Taylor expansions, we show that we actually solve at with a centered scheme in space the solution to:
\begin{equation}
  \begin{aligned}
    &\dxb \left( \rhoeb \ueb \right)-\Dx \dxb \left( \lbdmaxelec \dxb \rhoeb \right) = \Damk \left( \rhoeb - \Dt \dxb \left( \rhoeb \ueb \right) -\Dt\Dx \dxb \left( \lbdmaxelec \dxb \rhoeb \right)\right)\\
    \Leftrightarrow& \left( 1-\Dt\Damk \right)\dxb \left( \rhoeb \ueb \right)-\Dx \left( 1-\Dt\Damk \right)\dxb \left( \lbdmaxelec \dxb \rhoeb \right) = \Damk \rhoeb
  \end{aligned}
\end{equation}
In practice, we have seen in our simulations that $\Damk$ in the converged state was of the order of unity while typically $\Dt \approx 10^{-5}$, so that we have $\Damk\Dt\ll 1$ finally leading to the equations:
\begin{equation}\label{eq:Modifiedstst}
  \dxb \left( \rhoeb \ueb \right)-\Dx \dxb \left( \lbdmaxelec \dxb \rhoeb \right) = \Damk\rhoeb.
\end{equation}
We highlight that, if we do not neglect the terms factored by $\Dt\Damk$, the modified equation would depend on the operator splitting that we have chosen, an aspect of the numerical scheme that purely affects the time resolution of the method and that might have not been expected to have an influence on the converged and time constant solution of the numerical scheme.

The modified equation, i.e., \cref*{eq:Modifiedstst}, differs from the theoretical stationary form of \cref*{eq:ScldMassElec} because of the term $\Dx \dxb \left( \lbdmaxelec \dxb \rhoeb \right)$. Because of its form and its origin, this term is commonly called numerical diffusion. Theoretically, this term goes to zero as $\Dx\to0$. Nonetheless, there is in practice a limit to the smallest mesh size $\Dx$ that we can achieve, if only because we have to solve Poisson equation at each time step and that the numerical scheme quickly becomes more and more inefficient and computationally costly as the number of cell $N$ increases. We now try to estimate the weight of this additional term in the sheath area. Integrating \cref*{eq:Modifiedstst} from $0$ to a position $x$, and using the fact that both $\dxb \rhoeb$ and $\ueb$ are null at the center of the domain by symmetry, we obtain:
\begin{equation}\label{eq:ModSol}
  \rhoeb \ueb = \Dx \lbdmaxelec \dxb \rhoeb + \Damk \int_0^x \rhoeb \text{d}\xprime
\end{equation} 
instead of the actual solution of the steady state $\rhoeb \ueb = \Damk \int_0^x \rhoeb \text{d}\xprime$. Additionally, we use the approximation $\dxb \rhoeb \approx \rhoeb \dxb \phib$ and $\dxb \phib \approx \Vsheath /\debyeb$, leading to numerical error of the form
\begin{equation}\label{eq:NumDiff}
  err_{num} \propto \frac{\rhoeb\Dx}{\sqrt{\varepsilon\lambdaSq}}.
\end{equation}
This formula matches with our observation the numerical error increases for increasing atomic masses, i.e., increasing ion-to-electron mass ratios. It is also observed in numerical experiment that are not presented in this paper, for the sake of conciseness, that the numerical diffusion increases for increasing domain lengths  $\Lo$ as compared to $\debye$. This formula also explains why the numerical diffusion seems to peak approximately at the middle of the sheath and then is reduced in the vicinity of the wall, as it is proportional to the density of electron $\rhoeb$, which collapses at the boundaries. The phenomenon does not entirely vanishes at the wall, for two reasons: (1) $\rhoeb$ does not exactly vanish either at the wall, and (2) the discrete boundary conditions are not consistent with our problem, as identified above as problem 3.

In any case, if we wish to remove the non-physical peak of numerical diffusion in the sheath, we know from \cref*{eq:NumDiff} that in the test case that we consider we should have $\Dx \ll 0.02/( \sqrt{1836} \times 0.2 ) \approx 10^{-3}$, which would be extremely costly, even in 1D, since we have to solve Poisson at each time step, and even more unrealistic for 2D applications. We insist on the fact that these values are obtained for the lightest of gases, and would be amplified for heavier ones such as helium, argon or xenon, or for larger domains. Therefore, even going for second-order methods will not bring a satisfactory solution to problem 1.
\subsubsection*{Stability of the HLL solver for ions, influence of the source term}
We now move on to the analysis of problem 2, that is the lack of stability when the HLL approximate solver is used with ions at low values of $\kappa$. For the HLL solver to be stable, it is necessary (see for instance \cite{toro2013riemann,leveque2002finite}) that $\HLLbpion$ (and respectively $\HLLbmion$) be an upper bound (respectively a lower bound) to the velocities of propagation of information in the hyperbolic problem, meaning the velocities of propagation of acoustic waves for instance. Our choice of $\HLLbpion$ and $\HLLbmion$ does indeed obey this restriction for the convective part of the problem. However, in the present case, we know that ions are not an isolated fluid, but are coupled with electrons via electrostatic interactions. Information spreads at a much higher speed within the electron fluid. The result is that, as it is well known in the literature \cite{Chabert}, ion acoustic waves actually propagate typically at the Bohm velocity which much higher than their thermal velocity ($\vthib = \sqrt{\kappa} \ll 1 = \uBohmb$).

Indeed, within the bulk, where we can assume that the limit $\lambdaSq\ll1$ implies quasi-neutrality, that is $\rhoeb \approx \rhoionb$, and where we can neglect the electron macroscopic velocity $\ueb \ll \vtheb = \varepsilon^{-1/2}$, leading to the Boltzmann relation $\dxb \rhoeb = \rhoeb \dxb \phib$, we can replace the Lorentz force term for ions using the follwing approximations:
\begin{equation*}
  -\rhoionb\dxb\phib \approx -\rhoionb \frac{\dxb \rhoeb}{\rhoeb} \approx -\dxb \rhoionb
\end{equation*}
which when injected into the momentum equation of ions \cref*{eq:ScldMomIon} leads to :
\begin{equation*}
  \dtb \left( \rhoionb \uib \right) + \dxb \left( \rhoionb \uib^{2}+\rhoionb\left( \kappa + 1 \right) \right) = 0
\end{equation*}
so that we obtain a effective sound velocity for ions $\vthibeff = \sqrt{\kappa+1} = \sqrt{\vthib^{2}+\uBohmb^{2}}$.

Consequently, the HLL solver switches to its supersonic form, where it is essentially an upwind scheme, before the ions actually have a macroscopic velocity that is higher than the effective sound velocity in the plasma, that is approximately the Bohm speed $\uBohmb$ when $\kappa\ll1$, likely causing the instability observed in \cref*{fig:1stTCe1Rusi1HLLLie1}. The fact that instabilities appears in the bulk of the plasma, where the coupling of the two fluids is at its highest, and not in the sheath, where the two fluids are practically decoupled, underpins this interpretation of the phenomenon.
\subsubsection*{Discrete boundary conditions for electrons}
For the last problem, problem 3, we believe that the source of inaccuracies is the fact that the boundary conditions were derived without taking into account that the problem is not only a convective problem, but also features a source term. Similar loss of accuracy at the boundaries when an operator splitting is used has also been observed for instance for advection-diffusion problems, see \cite{hundsdorfer2003numerical} for a detailed analysis. In our case, we believe the discrete boundary conditions should include the effect of the Lorentz force to yield more accurate results.
\section{A novel accurate and stable first-order scheme}
In this section we gradually present how to correct each of the identified problems. We incrementally incorporate these variation of the original classical method to build a new method, and all numerical fixes are collected in our novel approach at the end of the section.

\subsection*{Adequate scaling of ion numerical diffusion for the HLL scheme in the quasi-neutral bulk}
\begin{figure}[htbp]
  \centering
  \includegraphics[width=1.\columnwidth]{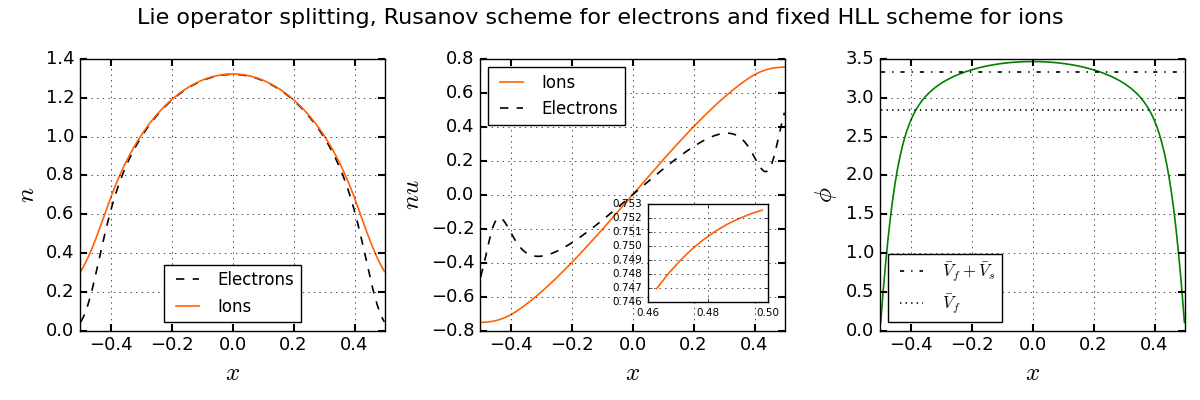}
  \caption{Number density, momentum and electric potential profiles at final time $\timef = 4\timeo$. Plasma parameters: $\varepsilon = 1/1836$, $\kappa = 0.0025$, $\debyeb = 0.02$. Numerical scheme: Lie splitting, Rusanov solver for electron and fixed HLL solver for ions and electrons. Simulation parameters: $N = 256$ cells, $\Dx \approx 4\times 10^{-3}$, $\Dt \approx 2\times 10^{-5}$.}
  \label{fig:1stTCe1Rusi1HLLfLie1}
\end{figure}

The first numerical fix that we present corresponds to problem 2, that is ion instabilities with the HLL solver in the quasi-neutral bulk of the plasma. As we explained in the previous section, the velocity of sound for the ion fluid is the Bohm velocity $\uBohmb$ and not their thermal velocity $\vthib = \sqrt{\kappa}$. Consequently, a simple fix to retrieve the stability of the solver consists in replacing $\vthib$ by the Bohm in the numerical diffusion of the scheme, that is setting $\HLLbpmion = \uRoeionb\pm \uBohmb$. The result obtained via this method are displayed in \cref*{fig:1stTCe1Rusi1HLLfLie1}.
As it can be seen in the figure, the numerical oscillations as well as the gradient at the last cell in the ion flux disappear with this new numerical scheme. In the following we refer to this scheme to as "fixed HLL".

\subsection*{Consistent discrete boundary conditions for electrons}
Secondly, we deal with problem 3, i.e., the problem of the consistency of the discrete boundary conditions that we employ for electrons. The condition that is imposed on the velocity of electrons is a direct consequence of the continuous problem boundary conditions and cannot be modified if we are to remain consistent with our model. Instead we focus on the boundary condition for electron density. To derive our new formula we first rewrite the momentum equation, \cref*{eq:ScldMomElec} in a non-conservative form:
\begin{equation}\label{eq:NCScldMomElec}
  \dtb \ueb+\dxb \left(\frac{\ueb^2}{2}+\frac{\ln\left(\rhoeb\right)}{\varepsilon}\right)=\frac{\dxb \phib}{\varepsilon}-\left(\nueb+\Damk\right)\ueb.
\end{equation}
By using that the fact that the desired boundary condition implies the condition $\dtb \ueb = 0$ on the walls, we obtain the discrete equation:
$$
\frac{\ln\left(\nelecGhost\right)-\ln\left(\nelecBndry\right)}{\varepsilon\Dx}= - \frac{\uelecGhost^2-\uelecBndry^2}{2\Dx}+\frac{\phiGhost-\phiBndry}{\varepsilon\Dx}-\left(\nueb+\Damk\right)\uwall
$$
where the index $G$ denotes the ghost cell and $B$ the cell at the boundary of the domain. This leads to the full set of equations:
\begin{subequations}\label{eq:HypBC}
    \begin{align}
        \nelecGhost &= \nelecBndry \exp \left( -\frac{\varepsilon}{2}\left(\uelecGhost^2-\uelecBndry^2\right)+\left(\phiGhost - \phiBndry\right)-\Dx\varepsilon\left(\nueb+\Damk\right)\uwall \right),\\
        \frac{\phiGhost+\phiBndry}{2} &=0,\\
        \frac{\uelecGhost+\uelecBndry}{2} &=\OutNorm \uwall.
    \end{align}
\end{subequations}
where $\OutNorm = 1$ for the right boundary and $\OutNorm = -1$ for the left boundary of the domain. As comparison, the previous boundary condition only used the formula $\nelecGhost = \nelecBndry$, leading to a slightly excessive value of the electron density at the wall. In turn this effect was amplified by the strong gradient in electron velocity near the boundaries, ultimately leading to inaccuracies in the electron momentum.

\begin{figure}[htbp]
  \centering
  \includegraphics[width=1.\columnwidth]{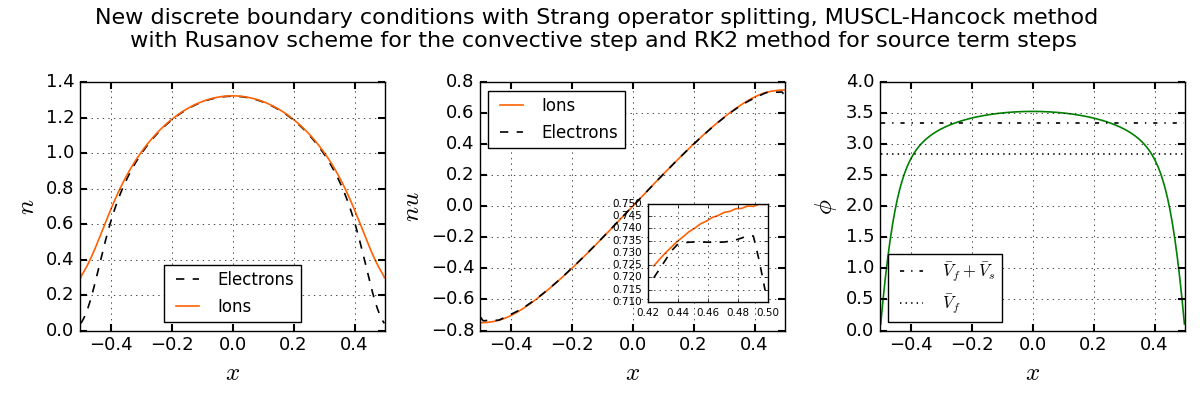}.
  \caption{Number density, momentum and electric potential profiles at final time $\timef = 4\timeo$. Plasma parameters: $\varepsilon = 1/1836$, $\kappa = 0.0025$, $\debyeb = 0.02$. Numerical scheme: new boundary conditions form \cref*{eq:HypBC}, Strang splitting, RK2 for the source term steps, MUSCL-Hancock methods with minmod slope limiters for the convective step and Rusanov solver for both electrons and first order scheme with Rusanov solver for ions. Simulation parameters: $N = 256$ cells, $\Dx \approx 4\times 10^{-3}$, $\Dt \approx 2\times 10^{-5}$.}
  \label{fig:1stTCe2Rusi1RusStrangBChyp}
\end{figure}

We display in \cref*{fig:1stTCe2Rusi1RusStrangBChyp} the results using this new boundary conditions with the Strang operator splitting, MUSCL-Hancock for the convective step with Rusanov solver for electrons and fixed HLL solver for ions and RK2 methods for the source term steps.
There is still a drop of momentum in the boundary cells but it has been reduced nearly by a factor eight. The solution is not perfectly smooth because although our boundary conditions are now consistent with the continuous problem, they do not exactly match the structure of the Strang operator splitting. Moreover, now that most of the numerical error stemming from the discretization of the boundary conditions have been significantly reduced, it is clear that going for second order reconstruction does not remove, in a satisfactory manner at least, the excessive diffusion of the scheme in the sheath area.

\subsection*{Controlled numerical diffusion for electron flux in the sheath}

The last step in devising our new scheme is to tackle problem 1, the most complex problem of the three, which is precisely the excessive numerical diffusion that appears in the sheath both with the first scheme and, to a lesser but still problematic extent, with the second order reconstruction for electron.

The method we propose will be essentially a small alteration of the Lie splitting presented in \cref*{eq:Lie1splitting}. The computation of the momenta and of the electric potential remains unchanged, only density profiles are obtained in a different way. In order to do so, the two first steps, \cref*{eq:Lie1ConvStep,eq:Lie1PoissonEq}, are not modified. Then, the computation of the ionization rate is delayed to first incorporate source terms for momentum equations. At that step, momentum values at time $\timenp$ are used to update number densities using a convective step with rescaled (at least for electrons) numerical diffusion. Lastly, the ionization rate is computed to compensate for the actual loss of ion mass in the domain and the ionization source term is added to the number densities of both ions and electrons. Overall, the unmodified part of method reads: 
\begin{subequations}\label{eq:Lie3splittingFirstPart}
  \begin{align}
    &\text{{\bfseries Step 1:} Convective step: }&\SVspeciesSj&=\SVspeciesnj-\frac{\Dt}{\Dx}\left(\DFlxspeciesnjph - \DFlxspeciesnjmh\right),\label{eq:Lie3ConvStep}\\
    &\text{{\bfseries Step 2:} Poisson equation: }&\phibVecnpj &= \PoissonOp \left(\nelecnionVecbSj\right),\label{eq:Lie3PoissonEq}\\
    &\text{{\bfseries Step 3:} Momentum source terms:}& \nUspeciesbnpj &= \nUspeciesbSj +\Dt \left( -\qspecies\frac{\phibnpjp-\phibnpjm}{2\Dx}-\nuspeciesbSj\nUspeciesbSj \right)
  \end{align}
\end{subequations}
Then we compute, using the updated values of momentums at time $\timenp$, the values of the number densities:  
\begin{subequations}\label{eq:Lie3splittingSecondPart}
  \begin{align}
    &\text{{\bfseries Step 4:} Particle flux convective step: }&\nspeciesbsSj&=\nspeciesbnj-\frac{\Dt}{\Dx}\left(\nUspeciesbnpjph - \nUspeciesbnpjmh\right),\label{eq:Lie3ConvStepPtwo}\\
    &\text{{\bfseries Step 5:} Ionization rate computation:}&\DamksS&=\IzOp\left( \nelecnionVecbsSnj \right),\label{eq:Lie3IzComp}\\
    &\text{{\bfseries Step 6:} Density source terms:}& \nspeciesbnpj &= \nspeciesbsSj + \Dt \DamksS \nelecbsSj,
  \end{align}
\end{subequations}
where the flux for the densities are computed using controlled diffusion, meaning that for electrons, we use the Rusanov flux where the numerical diffusion has been rescaled, following our analysis of the problem, leading to the new formula:
\begin{equation}
  \lbdmaxelec = \sqrt{\varepsilon\chi}\left( \max\left( \left|\uelecLb\right|, \left|\uelecRb\right| \right)+\vtheb \right)
\end{equation}
The scheme remains stable although we have severely reduced the numerical dissipation because using the flux at time $\timenp$ already brings some level of numerical diffusion. Similarly, we also scale, albeit in a little less dramatic fashion, the numerical diffusion of ions with the HLL solver to take into account that the scheme becomes naturally more diffusive when collisions are high. This leads to the modified coefficients:
\begin{equation}
  \HLLbpmion = \uRoeionb\pm \max\left( \Mscalingion\uBohmb, \vthib \right),
\end{equation}
where $\Mscalingion = 1/\left( 1+30\Dx\nuionb \right)$ is a scaling factor that was empirically tuned to reduce numerical viscosity in highly collisional regime when it would become excessive.

\begin{figure}[htbp]
  \centering
  \includegraphics[width=1.\columnwidth]{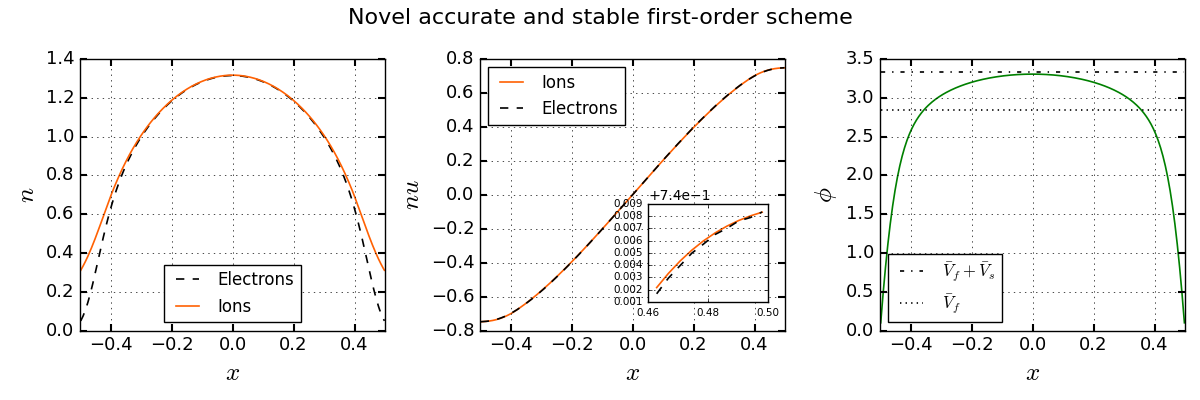}
  \caption{Number density, momentum and electric potential profiles at final time $\timef = 4\timeo$. Plasma parameters: $\varepsilon = 1/1836$, $\kappa = 0.0025$, $\debyeb = 0.02$. Numerical scheme: Modified Lie splitting, controlled diffusion Rusanov solver for electron and scaled diffusion fixed HLL solver for ions. Simulation parameters: $N = 256$ cells, $\Dx \approx 4\times 10^{-3}$, $\Dt \approx 2\times 10^{-5}$.}
  \label{fig:1stTCe1Rusi1HLLfLie3}
\end{figure}

The results of this new scheme are presented in \cref*{fig:1stTCe1Rusi1HLLfLie3}, the new boundary conditions from \cref*{eq:HypBC} were used.
We can observe that with this modified Lie splitting and controlled diffusion, the fluxes of electrons and ions are equal in the whole domain at the expected level of accuracy. Similarly, the electric potential reaches its theoretical value. Consequently, with minimal modifications of our original first-order scheme, we have obtained a method that outperforms classical second-order approaches available in the literature.


\section{Impact of the collisional regimes, gas composition, and spatial resolution on the novel scheme}
We now study the impact of the collisions, the electron-to-ion mass ratio and the spatial resolution on the proposed scheme.

\begin{figure}[htbp]
  \centering
  \includegraphics[width=1.\columnwidth]{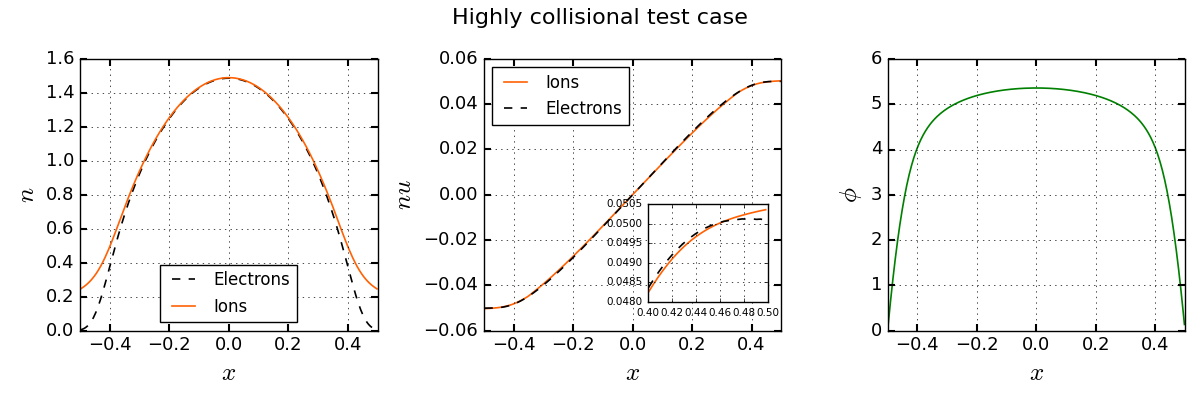}
  \caption{Number density, momentum and electric potential profiles at final time $\timef = 4\timeo$. Plasma parameters: $\varepsilon = 1/1836$, $\kappa = 0.0025$, $\debyeb = 0.02$, $\Lo/\lambdaion=10^{3}$, $\sigmane/\sigmanion = 10^{-1}$. Numerical scheme: Modified Lie splitting, controlled diffusion Rusanov solver for electron and scaled diffusion fixed HLL solver for ions. Simulation parameters: $N = 256$ cells, $\Dx \approx 4\times 10^{-3}$, $\Dt \approx 2\times 10^{-5}$.}
  \label{fig:2ndTCe1Rusi1HLLfLie3BChyp}
\end{figure}

First of all, we include the collisional terms in our model \cref*{eq:ScldEulerPoisson}. Our second test case it therefore to consider a highly collisional test case, where we choose a macro to mean path ratio for ions that is $\Lo/\lambdaion=10^{3}$ and as reference value for the ratio of cross sections we take $\sigmane/\sigmanion = 10^{-1}$. Result can be seen in \cref*{fig:2ndTCe1Rusi1HLLfLie3BChyp}.
We can see that the method still preserves the ambipolarity (equal ion and electron fluxes) to a very satisfactory level, and that the electric potential is correctly solved. Moreover, we also retrieve the characteristic sine profile for the densities in the bulk, which is expected from the literature \cite{Chabert}.

Secondly, we consider a heavier gas. For that purpose, we consider an argon gas, i.e. $\varepsilon = 1.36\times 10^{-5}$. The results are displayed in
\begin{figure}[htbp]
  \centering
  \includegraphics[width=1.\columnwidth]{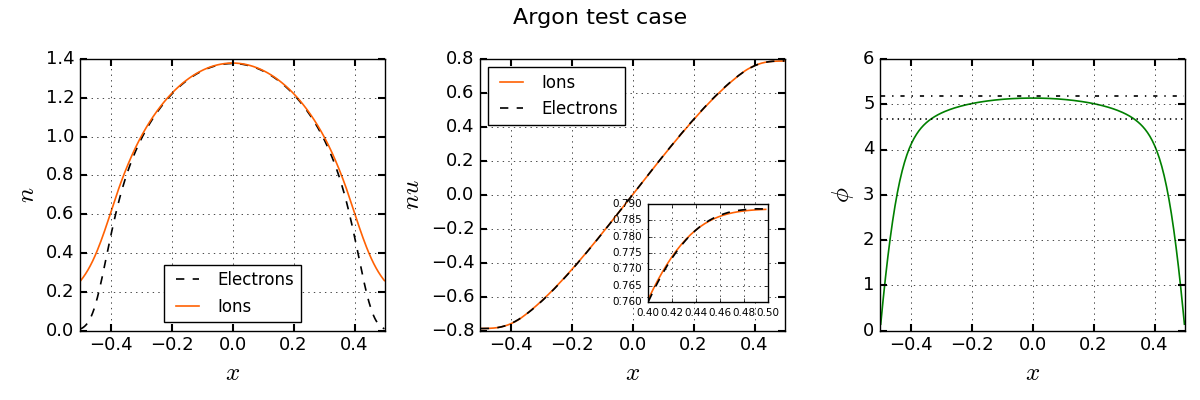}
  \caption{Number density, momentum and electric potential profiles at final time $\timef = 4\timeo$. Plasma parameters: $\varepsilon = 1.36\times 10^{-5}$, $\kappa = 0.0025$, $\debyeb = 0.02$, no collision. Numerical scheme: Modified Lie splitting, controlled diffusion Rusanov solver for electron and scaled diffusion fixed HLL solver for ions. Simulation parameters: $N = 256$ cells, $\Dx \approx 3\times 10^{-3}$, $\Dt \approx 4\times 10^{-6}$.}
  \label{fig:3ndTCe1Rusi1HLLfLie3BChyp}
\end{figure}
As it can be seen in the figure, the numerical scheme accurately represents the plasma-sheath transition even in the case of very small electron-to-ion mass ratios.

\begin{figure}[htbp]
  \centering
  \includegraphics[width=1.\columnwidth]{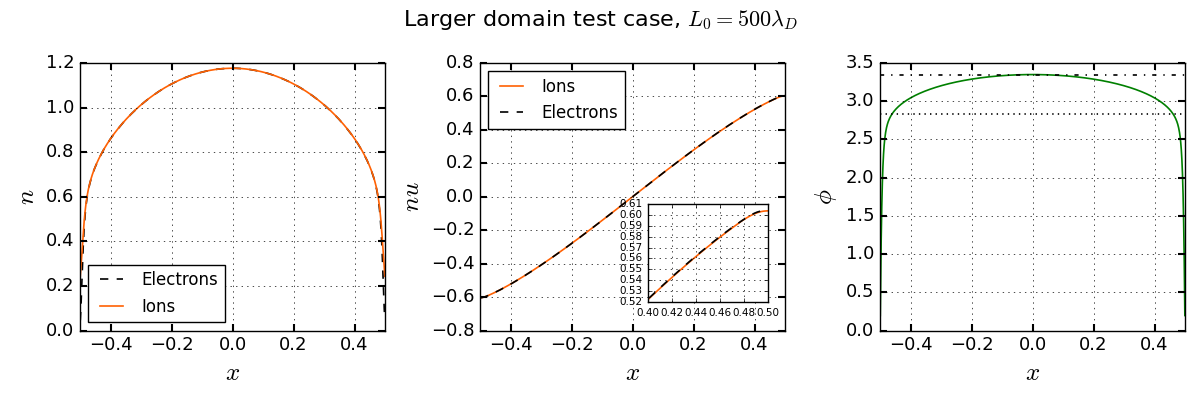}
  \caption{Number density, momentum and electric potential profiles at final time $\timef = 4\timeo$. Plasma parameters: $\varepsilon = 1/1836$, $\kappa = 0.0025$, $\debyeb = 0.002$, no collision. Numerical scheme: Modified Lie splitting, controlled diffusion Rusanov solver for electron and scaled diffusion fixed HLL solver for ions. Simulation parameters: $N = 1024$ cells, $\Dx \approx 1\times 10^{-3}$, $\Dt \approx 6\times 10^{-6}$.}
  \label{fig:4thTCe1Rusi1HLLfLie3BChyp}
\end{figure}

Lastly, we present a test case with a much larger domain with $\Lo = 500\debye$. The results are presented in \cref*{fig:4thTCe1Rusi1HLLfLie3BChyp}.
These results demonstrate that the method also behaves well when larger domains are considered, so long as we use at least two cells by Debye length.

\begin{figure}[htbp]
  \centering
  \includegraphics[width=1.\columnwidth]{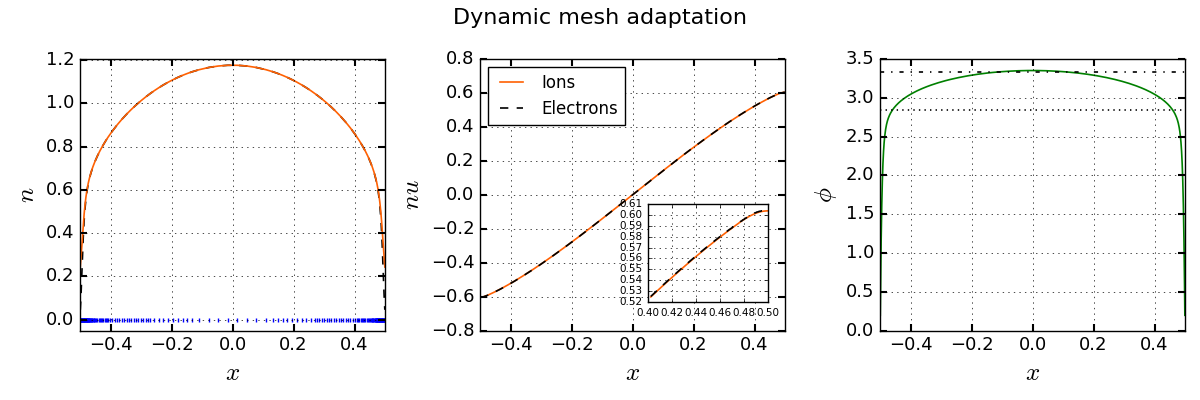}
  \caption{Number density, momentum and electric potential profiles at final time $\timef = 4\timeo$. Plasma parameters: $\varepsilon = 1/1836$, $\kappa = 0.0025$, $\debyeb = 0.002$, no collision. Numerical scheme: Modified Lie splitting, controlled diffusion Rusanov solver for electron and scaled diffusion fixed HLL solver for ions. Simulation parameters: $\Dx$ ranges from $2^{-5}$ to $2^{-10}$, $\Dt \approx 6\times 10^{-6}$. The vertical blue markers materialize the centers of the cells.}
  \label{fig:4thTCe1Rusi1HLLfLie3BChypMR}
\end{figure}

As a perspective for our work, we also run this last case with an adaptive mesh refinement method based on multi-resolution and error control \cite{DUARTE2015129,duarte:hal-00573043}, see \cref*{fig:4thTCe1Rusi1HLLfLie3BChypMR}.
We point out several facts. First of all, the method remains stable although the Debye length is not resolved in the center of the simulation, where $\Dx = 2^{-5}\geq 10\debyeb$. Second, the ambipolarity of the currents is maintained well on the whole domain and in particular in the sheath and does not seem to be affected by the fact that the mesh is no longer uniform. Third and last, the electric potential is well resolved and reaches the theoretical maximum value, which further validate our approach and paves the way for an efficient multi-dimension simulation.

\section{Conclusion}
In conclusion, we have proposed a minor modification of the classical Lie splitting that employed rescaled diffusion techniques to obtain a numerical method that performs well in various regimes to simulate bounded plasmas, significantly outperforming even second-order methods. The method could now be extended to a wide range of 1D applications and we have also been able to extend the scheme for 2D configuration, which we will present in future work.

\begin{acknowledgments}
We acknowledge the precious help of Lo\"ic Gouarin and the use of the samurai code he develops  (\url{https://github.com/hpc-maths/samurai}) within the framework of Initiative HPC@Maths (PI M. Massot and L. Gouarin - \url{https://initiative-hpc-maths.gitlab.labos.polytechnique.fr/site/}).

The first author of this proceeging is funded by a joint PhD grant from the French Ministry of Defence (AID - Agence Innovation D\'efense) and the R\'egion \^Ile-de-France (DIM MathInnov).
\end{acknowledgments}

\bibliography{Reboul_Alvarez_Laguna_Massot_1d_sheath_2022}
\appendix
\section{Strang splitting for the Euler-Poisson equations}
We describe here the Strang splitting \cite{toro2013riemann,leveque2002finite} that is used for our second order scheme. We assume that at the beginning of the step we have at our disposal $\phibVecnj$, that is the electric potential at time $\timen$. The structure of the scheme reads:
\begin{subequations}
  \begin{align}
  &\text{{\bfseries Step 1:} Source term half step, RK2:}&\Damkn &= \IzOpone\left( \nUiVecbnoneN, \nelecVecbnj \right)\\
  &&\SVspeciesnonej &= \SVspeciesnj+\frac{\Dt}{4} \St\left(\SVspeciesnj, \Damkn, \nelecbnj, \phibnjpm\right)\\
  &&\Damknone &= \IzOpone\left( \nUiVecbnonejoneN, \nelecVecbnonej \right)\\
  &&\SVspeciesntwoj &= \SVspeciesnj+\frac{\Dt}{2}\St\left(\SVspeciesnonej, \Damknone, \nelecbnonej, \phibnjpm\right)\\
  &\text{{\bfseries Step 2:} Convective step, MUSCL-Hancock:}&\SVspeciesnthreej &= \SVspeciesntwoj-\frac{\Dt}{\Dx}\left(\bMFlxspeciesntwojph-\bMFlxspeciesntwojmh\right)\\
  &\text{{\bfseries Step 3:} Poisson equation:}&\phibVecnpj &= \PoissonOp \left(\nelecnionVecbnthreej\right),\\
  &\text{{\bfseries Step 4:} Source term half step, RK2:}&\Damknthree &= \IzOpone\left( \nUiVecbnthreejoneN, \nelecVecbnthreej \right)\\
  &&\SVspeciesnfourj &= \SVspeciesnthreej+\frac{\Dt}{4} \St\left(\SVspeciesnthreej, \Damknthree, \nelecbnthreej, \phibnpjpm\right)\\
  &&\Damknfour &= \IzOpone\left( \nUiVecbnfourjoneN, \nelecVecbnfourj \right)\\
  &&\SVspeciesnpj &= \SVspeciesnthreej+\frac{\Dt}{2}\St\left(\SVspeciesnfourj, \Damknfour, \nelecbnfourj, \phibnpjpm\right)
\end{align}
\end{subequations}
\end{document}